\documentclass[12pt]{article}

\usepackage{amsmath,amsthm,amsfonts,amssymb,latexsym,amscd}

\font\fiverm=cmr5

 \input{prepictex.tex}
 \input{pictex.tex}
 \input{postpictex.tex}

\begin{document}

\newtheorem*{theo}{Theorem}
\newtheorem*{pro}{Proposition}
\newtheorem*{cor}{Corollary}
\newtheorem*{lem}{Lemma}
\newtheorem{theorem}{Theorem}[section]
\newtheorem{corollary}[theorem]{Corollary}
\newtheorem{lemma}[theorem]{Lemma}
\newtheorem{proposition}[theorem]{Proposition}
\newtheorem{conjecture}[theorem]{Conjecture}
\newtheorem{definition}[theorem]{Definition}
\newtheorem{problem}[theorem]{Problem}
\newtheorem{remark}[theorem]{Remark}
\newtheorem{example}[theorem]{Example}
\newcommand{\Naturali}{{\mathbb{N}}}
\newcommand{\Reali}{{\mathbb{R}}}
\newcommand{\Complessi}{{\mathbb{C}}}
\newcommand{\Toro}{{\mathbb{T}}}
\newcommand{\Relativi}{{\mathbb{Z}}}
\newcommand{\HH}{\mathfrak H}
\newcommand{\KK}{\mathfrak K}
\newcommand{\LL}{\mathfrak L}
\newcommand{\Cs}{$C^*$-algebra}
\newcommand{\as}{\ast_{\sigma}}
\newcommand{\tn}{\vert\hspace{-.3mm}\vert\hspace{-.3mm}\vert}
\def\mA{{\mathfrak A}}
\def\A{{\mathcal A}}
\def\mB{{\mathfrak B}}
\def\B{{\mathcal B}}
\def\C{{\mathcal C}}
\def\D{{\mathcal D}}
\def\F{{\mathcal F}}
\def\H{{\mathcal H}}
\def\J{{\mathcal J}}
\def\K{{\mathcal K}}
\def\L{{\cal L}}
\def\N{{\cal N}}
\def\M{{\cal M}}
\def\O{{\mathcal O}}
\def\P{{\cal P}}
\def\S{{\cal S}}
\def\T{{\cal T}}
\def\U{{\cal U}}
\def\W{{\cal W}}
\def\b{\lambda_B(P}
\def\j{\lambda_J(P}
\def\Z{{\mathbb{Z}}}

\title{Endomorphisms of ${\mathcal O}_n$ which preserve the canonical
  UHF-subalgebra}

\author{Roberto Conti, Mikael R{\o}rdam, Wojciech Szyma{\'n}ski}

\date{5 October 2009}
\maketitle

\renewcommand{\sectionmark}[1]{}

\vspace{7mm}
\begin{abstract}
Unital endomorphisms of the Cuntz algebra $\O_n$ which preserve the canonical
UHF-subalgebra $\F_n \subseteq \O_n$ are investigated. We give examples of
such endomorphisms $\lambda = \lambda_u$ for which the associated
unitary element $u$ in $\O_n$ (which satisfies $\lambda(S_j) = uS_j$
for all $j$) does not belong to $\F_n$. One such example, in the case
where $n=2$,  arises from a construction of a unital endomorphism  on $\O_2$
which preserves the canonical UHF-subalgebra and where the relative commutant of its
image in $\O_2$ contains a copy of $\O_2$.
\end{abstract}

\vfill

\vspace{3mm}

\newpage

\section{Introduction}
The study of endomorphisms of Cuntz algebras continues to attract attention
of researchers. On the one hand, such endomorphisms naturally arise in a
number of contexts including
index theory and subfactors, entropy, and classical dynamical systems
on the Cantor set.
On the other hand, they exhibit interesting and intriguing features while being
concrete enough to allow explicit albeit sometimes binding computations.

It is a fundamental fact that there is a one-to-one correspondence
between unitaries in $\O_n$ and unital endomorphisms on $\O_n$ whereby
$u$ in $\O_n$ corresponds to the endomorphism $\lambda_u$ which maps
the $j$th canonical generator $S_j$ of $\O_n$ onto $uS_j$ for
$j=1,2,\dots, n$. In the ground breaking paper by Cuntz on this
subject, \cite{Cun2}, it is noted that $\lambda_u$ maps the canonical
UHF-subalgebra $\F_n$ of $\O_n$ into itself whenever $u$ belongs to
$\F_n$; and the question if the converse also holds is considered. This
indeed is true in many cases (as one can deduce from \cite{Cun2}), for
example  if one knows in advance that the range of the endomorphism $\lambda_u$
is globally invariant under the gauge action of $\mathbb T$. This assumption is
already sufficient to cover several interesting cases, e.g. if
$\lambda_u$ is an automorphism of $\O_n$.

We show in this paper that this converse statement is false
in general, i.e., there is a unitary element $u$ in $\O_n$ which does not
belong to $\F_n$ but where $\lambda_u$ maps $\F_n$ into itself.

The paper is organized in the following way.
In section \ref{counterexample}, after some preliminaries, we present
a general framework for
finding the announced counterexamples and we discuss a specific
example in the case of $\O_2$ that arises in a combinatorial way.
In section \ref{large}, we exhibit a unitary $u$ in the UHF-subalgebra of
$\O_2$ such that the image of the  corresponding endomorphism $\lambda_u$ has
relative commutant containing a copy of $\O_2$. One can then easily find another
unitary $v$ in $\O_2$ such that $v$ does not belong to $\F_2$ but where
$\lambda_v$ agrees with $\lambda_u$ on $\F_2$, whence in particular $\lambda_v$ maps
$\F_2$ into itself. From this construction one gets as a byproduct an embedding
of $\O_2 \otimes \O_2$ into $\O_2$ that maps $\F_2 \otimes \F_2$ into $\F_2$.

It remains an interesting open problem if for every unital
endomorphism $\lambda$ on $\O_n$ which maps $\F_n$ into itself there
exists a unitary element $u$ in $\F_n$ such that $\lambda$ and
$\lambda_u$ agree on $\F_n$.

In section \ref{epcs}, we expand our initial observations on
endomorphisms preserving the canonical
UHF-subalgebra in a more systematic manner.
In section \ref{sncocycles}, we study a particularly interesting class
of such endomorphisms related
to certain elements in the normalizer of the canonical MASA.

Finally, we would like to mention that endomorphisms preserving the core $AF$-subalgebras of
certain $C^*$-algebras corresponding to rank-2 graphs (generalizing the Cuntz algebras)
have been very recently considered in \cite{Y}.

\vspace{3mm}\noindent{\bf Acknowledgements.}
The need for clarification of some of the issues considered in the present paper was raised
in a discussion with Adam Skalski and Joachim Zacharias, who we warmly thank. The first mentioned
author is grateful to Uffe Haagerup for supporting his visit to Odense in June 2009, where part
of this work has been done.

\section{A counterexample}\label{counterexample}

If $n$ is an integer greater than 1, then the Cuntz algebra $\O_n$ is a unital, simple
$C^*$-algebra generated by $n$ isometries $S_1, \ldots, S_n$, satisfying
$\sum_{i=1}^n S_i S_i^* = I$, \cite{Cun1}.
We denote by $W_n^k$ the set of $k$-tuples $\alpha = (\alpha_1,\ldots,\alpha_k)$
with $\alpha_m \in \{1,\ldots,n\}$, and by $W_n$ the union $\cup_{k=0}^\infty W_n^k$,
where $W_n^0 = \{0\}$. We call elements of $W_n$ multi-indices.
If $\alpha = (\alpha_1,\ldots,\alpha_k) \in W_n$, then $S_\alpha = S_{\alpha_1} \ldots S_{\alpha_k}$
($S_0 = I$ by convention) and $P_\alpha=S_\alpha S_\alpha^*$.
Every word in $\{S_i, S_i^* \ | \ i = 1,\ldots,n\}$ can be uniquely expressed as
$S_\alpha S_\beta^*$, for $\alpha, \beta \in W_n$ \cite[Lemma 1.3]{Cun1}.
If $\alpha \in W_n^k$ then $|\alpha| = k$ is the length of $\alpha$.

$\F_n^k$ is the $C^*$-algebra generated by all words of the form
$S_\alpha S_\beta^*$, $\alpha, \beta \in W_n^k$, and it is isomorphic to the
matrix algebra $M_{n^k}({\mathbb C})$. $\F_n$, the norm closure of
$\cup_{k=0}^\infty \F_n^k$, is the UHF-algebra of type $n^\infty$,
called the core UHF-subalgebra of $\O_n$, \cite{Cun1}. It is the fixed point algebra
for the periodic gauge action of the reals: $\alpha:{\mathbb R}\rightarrow{\rm Aut}(\O_n)$ defined
on generators as $\alpha_t(S_i)=e^{it}S_i$, $t \in {\mathbb R}$.

We denote by $\S_n$ the group of those unitaries in $\O_n$ which can be written
as finite sums of words, i.e., in the form $u = \sum_{j=1}^m S_{\alpha_j}S_{\beta_j}^*$
for some $\alpha_j, \beta_j \in W_n$. It turns out that $\S_n$ is isomorphic to
the Higman-Thompson group $G_{n,1}$ \cite{Ne}. We also denote $\P_n=\S_n\cap\U(\F_n)$. Then
$\P_n=\cup_k\P_n^k$, where $\P_n^k$ are permutation unitaries in $\U(\F_n^k)$.
That is, for each $u\in\P_n^k$ there is a unique permutation $\sigma$ of multi-indices
$W_n^k$ such that $u = \sum_{\alpha \in W_n^k} S_{\sigma(\alpha)} S_\alpha^*$.

For $u$ a unitary in $\O_n$ we denote by $\lambda_u$ the unital endomorphism of $\O_n$
determined by $\lambda_u(S_i) = u S_i$, $i=1,\ldots, n$. We denote by $\varphi$ the
canonical shift: $\varphi(x)=\sum_i S_ixS_i^*$, $x\in\O_n$. Note that $\varphi$
commutes with the action $\alpha$. If $u\in\U(\O_n)$ then for each positive integer $k$ we denote
$$ u_k = u \varphi(u) \cdots \varphi^{k-1}(u). $$
We agree that $u_k^*$ stands for $(u_k)^*$. If
$\alpha$ and $\beta$ are multi-indices of length $k$ and $m$, respectively, then
$\lambda_u(S_\alpha S_\beta^*)=u_kS_\alpha S_\beta^*u_m^*$. This is established through
a repeated application of the identity $S_i a = \varphi(a)S_i$, valid for all
$i=1,\ldots,n$ and $a \in \O_n$.

\begin{proposition}\label{uv}
Let $u$ be a unitary in $\O_n$ and let $v$ be a unitary in the
relative commutant
$\lambda_u(\F_n)' \cap \O_n$. Define $w := u \varphi(v)$. Then the
restrictions of
endomorphisms $\lambda_u$ and $\lambda_w$ coincide on
$\F_n$. Likewise, if $\tilde{w}=vu$
then the restrictions of endomorphisms $\lambda_u$ and
$\lambda_{\tilde{w}}$ coincide on $\F_n$.
\end{proposition}
\begin{proof}
It is enough to compute the action of $\lambda_w$ on all elements of the form
$S_{\alpha_1} \ldots S_{\alpha_k} S_{\beta_k}^* \ldots S_{\beta_1}^*$ for
every integer $k \geq 1$ and all $\alpha_i$ and $\beta_j$ in $\{1,\ldots,n\}$
for all $1 \leq i,j \leq k$. To this end, we verify by induction on $k$ that
$$ \lambda_w(S_{\alpha_1} \ldots S_{\alpha_k} S_{\beta_k}^* \ldots
S_{\beta_1}^*) =
   \lambda_u(S_{\alpha_1} \ldots S_{\alpha_k} S_{\beta_k}^* \ldots
   S_{\beta_1}^*). $$
Indeed, for $k=1$ we have
$$ \lambda_w(S_{\alpha_1}S^*_{\beta_1})=wS_{\alpha_1}S^*_{\beta_1}w^*=
   u\varphi(v)S_{\alpha_1}S^*_{\beta_1}\varphi(v)^*u^*=
   uS_{\alpha_1}S^*_{\beta_1}u^*=\lambda_u(S_{\alpha_1}S^*_{\beta_1}), $$
since $\varphi(v)$ and $S_{\alpha_1}S^*_{\beta_1}$ commute. Now
assuming the identity
holds for $k-1$, we have
\begin{align*}
\lambda_w(S_{\alpha_1} \ldots S_{\alpha_k} S_{\beta_k}^* \ldots S_{\beta_1}^*) & =
\lambda_w(S_{\alpha_1})\lambda_w(S_{\alpha_2} \ldots S_{\alpha_k} S_{\beta_k}^*
\ldots S_{\beta_2}^*)\lambda_w(S_{\beta_1})^* \\
& = u\varphi(v)S_{\alpha_1}\lambda_u(S_{\alpha_2} \ldots S_{\alpha_k} S_{\beta_k}^*
\ldots S_{\beta_2}^*)S^*_{\beta_1}\varphi(v)^*u^* \\
& = uS_{\alpha_1}v\lambda_u(S_{\alpha_2} \ldots S_{\alpha_k} S_{\beta_k}^*
\ldots S_{\beta_2}^*)v^*S^*_{\beta_1}u^* \\
& = uS_{\alpha_1}\lambda_u(S_{\alpha_2} \ldots S_{\alpha_k} S_{\beta_k}^*
\ldots S_{\beta_2}^*)S^*_{\beta_1}u^* \\
& = \lambda_u(S_{\alpha_1} \ldots S_{\alpha_k} S_{\beta_k}^* \ldots S_{\beta_1}^*),
\end{align*}
since $v$ is in the commutant of $\lambda_u(\F_n)$. The proof of the remaining claim is
similar.
\end{proof}

\begin{corollary}\label{Countercuntz}
Under the hypothesis of Proposition \ref{uv}, assume further that $u\in\F_n$. Then
$\lambda_u(\F_n)\subseteq\F_n$ and thus $\lambda_w(\F_n)\subseteq\F_n$.
However, $w$ belongs to $\F_n$ if and only if $v$ does.
\end{corollary}

The crucial role in the above construction is played by $\lambda_u(\F_n)' \cap \O_n$.
It turns out that this relative commutant can be calculated as follows (compare
\cite[Proposition 3.1]{L}).

\begin{proposition} \label{Prc}
Let $u$ be a unitary in $\O_n$, then
\begin{equation}\label{rcUHF}
\lambda_u(\F_n)' \cap \O_n = \bigcap_{k \geq 1} ({\rm Ad} u\circ \varphi)^k (\O_n).
\end{equation}
\end{proposition}

\begin{proof}
Clearly an element $x \in \O_n$ lies in
$\lambda_u(\F_n)' \cap \O_n$ if and only if, for all $k \geq 1$ and all $y \in \F_n^k$,
$x$ commutes with $\lambda_u(y) = u_k y u_k^*$, i.e.
$$u_k^* x u_k \in (\F_n^k)' \cap \O_n = \varphi^k(\O_n) \ . $$
This means precisely that, for each $k \geq 1$,
$x$ lies in the range of ${\rm Ad}(u_k) \varphi^k = ({\rm Ad}u \circ \varphi)^k$.
\end{proof}

It is also useful to observe that
${\rm Ad} u\circ \varphi$ restricts to an automorphism of $\lambda_u(\F_n)' \cap \O_n$.
This follows from the following simple lemma.

\begin{lemma}\label{Laut}
Let $\mathfrak A$ be a unital $C^*$-algebra and $\rho$ an injective unital $*$-endomorphism
of $\mathfrak A$, then $\rho$ restricts to a $*$-automorphism of
$${\mathfrak A}_\rho := \bigcap_{k \in {\mathbb N}} \rho^k({\mathfrak A}) \ . $$
\end{lemma}
\begin{proof}
One has a descending tower of unital $C^*$-subalgebras of $\mathfrak A$,
$${\mathfrak A} \supset \rho({\mathfrak A}) \supset \rho^2({\mathfrak A}) \supset \ldots \ , $$
thus ${\mathfrak A}_\rho$ is a unital $C^*$-subalgebra of $\mathfrak A$.
An element $x \in {\mathfrak A}_\rho$ satisfies
$$x = \rho(x_1) = \rho^2(x_2) = \cdots = \rho^k(x_k) = \ldots$$
for elements $x_1, \ldots, x_k, \ldots$ in $\mathfrak A$.
It is then clear that $\rho$ maps ${\mathfrak A}_\rho$ into itself,
and moreover $x_1, \ldots, x_k, \ldots \in {\mathfrak A}_\rho$ so that in particular
$\rho({\mathfrak A}_\rho) = {\mathfrak A}_\rho$.
\end{proof}
Endomorphisms $\rho$ for which ${\mathfrak A}_\rho = {\mathbb C}1$ are often called {\it shifts}.

\medskip

Corollary \ref{Countercuntz} shows how to construct examples of unitaries 
$w$ outside $\F_n$
for which nevertheless $\lambda_w(\F_n) \subseteq \F_n$. To this end, it suffices to find a unitary
$u\in\F_n$ such that the relative commutant $\lambda_u(\F_n)' \cap \O_n$ is not contained in
$\F_n$. This is possible. In fact, one can even find unitaries in a matrix algebra $\F_n^k$
such that $\lambda_u(\O_n)'\cap\O_n$ is not contained in $\F_n$. The existence of such unitaries
was demonstrated in \cite{CP}. The relative commutant $\lambda_u(\O_n)'\cap\O_n$
coincides with the space $(\lambda_u,\lambda_u)$ of self-intertwiners of the
endomorphism $\lambda_u$, which can be computed as
$$ (\lambda_u,\lambda_u) = \{ x\in\O_n:x=({\rm Ad} u\circ\varphi)(x)\} \ . $$

\begin{example}\label{countercuntz}
\rm We give an explicit example of
a permutation unitary $u\in\P_2^4$ and a unitary $v$ in $\S_2\setminus\P_2$ such that
$v\in(\lambda_u,\lambda_u)$. Indeed, one can check by a lengthy but straightforward
computation that the pair:
\begin{align*}
u & = S_1 S_1 S_1 S_1^* S_1^* S_1^* + S_1 S_2 S_1 S_1 S_1^* S_2^* S_1^* S_1^* +
S_1 S_2 S_2 S_1 S_2^* S_2^* S_1^* S_1^*
\\
& + S_2 S_1 S_1 S_1 S_1^* S_1^* S_2^* S_1^* + S_1 S_2 S_1 S_2 S_2^* S_1^* S_2^* S_1^* +
S_2 S_2 S_1 S_2^* S_2^* S_1^* \\
& + S_1 S_1 S_2 S_1^* S_1^* S_2^* + S_2 S_1 S_2 S_1 S_1^* S_2^* S_1^* S_2^* +
S_1 S_2 S_2 S_2 S_2^* S_2^* S_1^* S_2^* \\
& + S_2 S_1 S_1 S_2 S_1^* S_1^* S_2^* S_2^* + S_2 S_1 S_2 S_2 S_2^* S_1^* S_2^* S_2^* +
S_2 S_2 S_2 S_2^* S_2^* S_2^*
\in \P_2^4
\end{align*}
\begin{align*}
v & = S_1 S_2 S_2 S_1^* S_1^* + S_1 S_1 S_1 S_1^* S_2^* S_1^* + S_2 S_1 S_1 S_2^* S_2^* S_1^* \\
& + S_2 S_2 S_1^* S_1^* S_2^* + S_1 S_1 S_2 S_2^* S_1^* S_2^* + S_1 S_2 S_1 S_1^* S_2^* S_2^* \\
& + S_2 S_1 S_2 S_2^* S_2^* S_2^* \in \S_2 \setminus \P_2
\end{align*}
does the job.
\end{example}

The way the examples in this and the following section have been constructed leaves open
the possibility that for each endomorphism of $\O_n$ globally preserving the core
UHF-subalgebra there exists another one, induced by a unitary in $\F_n$, which restricts to
the same endomorphism of $\F_n$. So far, this question has not been settled in full
generality and we would like to leave it as an open problem.

\section{An endomorphism on $\O_2$ with relative commutant containing
$\O_2$} \label{large}

\noindent As observed in the previous section,
if $u$ is a unitary element in a
Cuntz algebra $\O_n$ and $v$ is a unitary element in the relative
commutant $\lambda_u(\F_n)' \cap \O_n$, then $\lambda_u$ and
$\lambda_{vu}$ will agree on the canonical UHF-algebra $\F_n$
contained in $\O_n$. Hence, if $u$ belongs to $\F_n$
and $\lambda_u(\F_n)' \cap \O_n$ is not contained in $\F_n$, then
one can choose $v$ as above such that $v$, and hence $vu$, do not
belong to $\F_n$; whereas $\lambda_{vu}$ will map $\F_n$ into
itself. We constructed an example of such a unitary element
$u$ in Example~\ref{countercuntz} above. In this section we shall
construct another example, in the case where $n=2$, where the relative
commutant $\lambda_u(\F_2)' \cap \O_2$ contains $\O_2$ and therefore
is not contained in $\F_2$.

It is well-known, \cite{Ror:o2o2}, that $\O_2$ is isomorphic to $\O_2
\otimes \O_2$. In particular there is a unital embedding $\O_2
\otimes \O_2 \to \O_2$. If one composes that with the embedding
$\O_2 \to \O_2 \otimes \O_2$ given by $x \mapsto x \otimes 1$, then
one obtains an endomorphism $\lambda \colon \O_2 \to \O_2$ such that
$\lambda(\O_2)' \cap \O_2$ contains a unital copy of $\O_2$. We
show that one can choose this endomorphism $\lambda$ such that it is
of the form $\lambda = \lambda_u$ for some unitary $u$ in $\F_2$.

Let $\eta \colon \O_2 \to D$ be a unital $*$-homomorphism. Let $S_1,S_2$ be the
two canonical generators of $\O_2$. Define unital endomorphisms $\varphi$
on $\O_2$ and $\psi$ on $D$ by
\begin{equation} \label{eq:mu}
\varphi(x) = S_1xS_1^* + S_2xS_2^*, \qquad \psi(y) =
\eta(S_1)y\eta(S_1)^* + \eta(S_2)y\eta(S_2)^*,
\end{equation}
for $x \in \O_2$ and $y \in D$. One has that $\psi \circ \eta =
\eta \circ \varphi$. Close inspection of the proof of Theorem~3.6 from
\cite{Ror:cuntz} shows that the following holds:

\begin{theorem}[cf.\ \cite{Ror:cuntz}] \label{thm:r}
Let $D$ be a
 unital properly infinite $C^*$-algebra, let $\eta \colon \O_2 \to D$ be a
 unital $*$-homomorphism, and let $D_0$ be a unital sub-$C^*$-algebra of $D$ such that
\begin{description}
\item[\hspace{2mm}{\rm (i)}] $D_0$ is $K_1$-injective and has bounded exponential length,
\item[\hspace{1mm}{\rm (ii)}] $D_0$ is invariant under the endomorphism $\psi$ on $D$
  associated with $\eta$ (as defined in \eqref{eq:mu} above).
\item[{\rm (iii)}] $D_0$ contains $\eta(\F_2)$.
\end{description}
It follows that $\{v\psi(v)^* \mid v \in \U(D_0)\}$ is dense in
$\U(D_0)$.
\end{theorem}

\noindent Combining the theorem above with \cite[Lemma 1]{Ror:o2o2} we get
the following:

\begin{proposition} \label{prop:1}
There is a sequence $\{v_n\}$ of unitaries in
  $\F_2$ such that the corresponding sequence $\{\lambda_{v_n}\}$ of
  endomorphisms on $\O_2$ is asymptotically central.
\end{proposition}

\begin{proof} As in the proof of  \cite[Lemma 1]{Ror:o2o2}, if $\{v_n\}$ is a
  sequence of unitaries in $\O_2$, then $\{\lambda_{v_n}\}$ is
  asymptotically central if and only if
$$\varphi(v_n)^*v_n \to \sum_{i,j=1}^2 S_iS_jS_i^*S_j^*.$$
The unitary on the right-hand side belongs to $\F_2$. The existence of
the desired sequence $\{v_n\}$ of unitaries in $\F_2$ therefore
follows from Theorem~\ref{thm:r} with $D=\O_2$, $\eta= {\rm Id}$, and with
$D_0 = \F_2$.
\end{proof}

\begin{proposition} \label{prop:1.5} In $\O_2 \otimes \O_2$ consider
  the unitary element
$$u_0 = S_1^* \otimes S_1 + S_2^* \otimes S_2,$$
and let $B$ be the $C^*$-algebra generated by $\F_2 \otimes \F_2 \cup
\{u_0\}$. Then
$$B = C^*(\F_2 \otimes \F_2,u_0) \cong
\big(\bigotimes_{n \in {\mathbb Z}} M_2\big) \rtimes_{\mathrm{shift}}
{\mathbb Z},$$
whence $B$ is a simple A$\mathbb T$-algebra of real rank zero. In
particular, $B$ is $K_1$-injective and has finite exponential rank.
Also, $B$ is
invariant under the endomorphism $\varphi \otimes {\rm Id}$ on $\O_2 \otimes
\O_2$ (where $\varphi$ is as defined in \eqref{eq:mu})

The ``half flip'' on $\O_2 \otimes \O_2$ is approximately
inner with unitaries belonging to $B$, i.e., there is a sequence
$\{z_n\}$ of unitaries in $B$ such that $z_n(x \otimes 1)z_n^* \to 1
\otimes x$ for all $x \in \O_2$.
\end{proposition}

\begin{proof} With $\varphi$ as above, put
$$E_{ij}^{(n)} = \begin{cases} \varphi^{(-n)}(S_iS_j^*) \otimes 1, &
n \le 0, \\ 1 \otimes \varphi^{n-1}(S_iS_j^*), & n \ge 1, \end{cases}
$$
for $i,j=1,2$ and for $n \in \Z$.
Then $C^*\big(E_{ij}^{(n)} \mid i,j=1,2\big)$, $n \in \Z$, is a com\-muting
family of \Cs  s each isomorphic to $M_2$,
and $$\F_2 \otimes \F_2 = C^*\big(E_{ij}^{(n)} \mid i,j=1,2, n
\in \Z\big).$$ Moreover, $u_0E_{ij}^{(n)}u_0^* = E_{ij}^{(n+1)}$
for all $i,j=1,2$ and for all $n \in \Z$.
This proves that
$$C^*(\F_2 \otimes \F_2,u_0) \cong \big(\bigotimes_{n \in {\mathbb Z}} M_2
\big) \rtimes_{\mathrm{shift}} {\mathbb Z}.$$
(Off hand, without going into issues of proper outerness of the shift
action, one would only get that the \Cs{} on the left-hand side is a
quotient of the crossed product on the right-hand side, but since the
latter is simple, we get the isomorphism for free.)

The crossed product on the right-hand side
is known to be a simple A$\mathbb T$-algebra of real rank
zero, \cite{BraKisRorSto:shift}. We must show that $C^*(\F_2 \otimes
\F_2,u_0)$ is invariant under $\kappa := \varphi \otimes {\rm
  Id}$. It is clear that $\F_2 \otimes
\F_2$ is invariant under $\kappa$.  A brief calculation shows that
$\kappa(u_0) = zu_0$, where
$$z = S_1S_1^* \otimes S_1S_1^* + S_1S_2^* \otimes S_2S_1^* + S_2S_1^*
\otimes S_1S_2^* + S_2S_2^* \otimes S_2S_2^*.$$
Notice that $z$ belongs to $\F_2 \otimes \F_2$. Therefore $\kappa(u_0)$
belongs to $C^*(\F_2 \otimes \F_2,u_0)$. As $\kappa^{n}(u_0) =
\kappa^{n-1}(z)\kappa^{n-1}(u_0)$ it follows by induction
that $\kappa^n(u_0)$ belongs to $C^*(\F_2 \otimes
\F_2,u_0)$ for all $n$. This proves that $C^*(\F_2 \otimes \F_2,u_0)$ is
invariant under $\kappa$.

It follows from Theorem~\ref{thm:r}, with $D =
\O_2 \otimes \O_2$, with $D_0 = B$, and with $\eta \colon \O_2
\to \O_2 \otimes \O_2$ given by $\eta(x) = x
\otimes 1$, that there is a sequence $\{z_n\}$ of unitaries in $B$
such that $z_n\kappa(z_n)^* \to u_0$. It is straightforward to check
that $1 \otimes S_j = u_0(S_j \otimes 1)$ and a standard calculation,
cf.\ \cite{Ror:cuntz}, then shows that $z_n(S_j \otimes 1)z_n^* \to 1
\otimes S_j$ for $j=1,2$.
\end{proof}

\begin{corollary} \label{cor:1.5} Let $D$ be a unital $C^*$-algebra, and suppose
    that $\eta_1, \eta_2 \colon
\O_2 \to D$ are unital $*$-homomorphisms with commuting
    images. There is a sequence $\{w_n\}$ of unitaries in
the sub-\Cs{} $D_0 = C^*(\eta_1(\F_2),u)$, where
$$u =
    \eta_2(S_1)\eta_1(S_1)^*+\eta_2(S_2)\eta_1(S_2)^*,$$
such that $w_n\eta_1(x)w_n^* \to
\eta_2(x)$ for all $x \in \O_2$.
\end{corollary}

\begin{proof} The $*$-homomorphisms $\eta_1$ and $\eta_2$ induce a
  $*$-homomorphism
  $\eta \colon \O_2 \otimes \O_2 \to D$ given by
$$\eta(x \otimes y) = \eta_1(x)\eta_2(y), \qquad x,y \in
\O_2.$$
In the notation of Proposition~\ref{prop:1.5} we have
$$\eta(u_0) = u, \qquad \eta(\F_2 \otimes 1) =
\eta_1(\F_2), \qquad  \eta(1 \otimes \F_2) =
\eta_2(\F_2).$$
It follows from Proposition~\ref{prop:1.5} and its proof that $1
\otimes \F_2$ is contained in the \Cs{} generated by $\{E_{ij}^{(0)}\}$
and $u_0$ and hence is contained
in $C^*(\F_2 \otimes 1, u_0)$. The \Cs{} $B$ from
that proposition is therefore generated by $\F_2 \otimes 1$ and $u_0$,
which shows that  $\eta(B) = D_0$.

Let $\{z_n\}$ be as in Proposition~\ref{prop:1.5} and put $w_n=
\eta(z_n) \in D_0$. Then
$$w_n \eta_1(x)w_n^* = \eta(z_n (x \otimes 1)z_n^*) \to
\eta(1 \otimes x) = \eta_2(x)$$
for all $x \in \O_2$.
\end{proof}

\begin{proposition} \label{prop:2}
There are sequences $\{v_n\}$ and $\{w_n\}$ of
  unitaries in $\F_2$ such that
\begin{description}
\item[\hspace{2mm}{\rm (i)}] $\{\lambda_{v_n}\}$ is asymptotically central in $\O_2$,
\item[\hspace{1mm}{\rm (ii)}] $\|w_n \lambda_{v_{n+1}}(S_j) w_n^* - \lambda_{v_n}(S_j)\| <
  2^{-n}$ for all $n \in {\mathbb N}$ and for $j=1,2$,
\item[{\rm (iii)}] $\|w_nS_jw_n^*-S_j\| < 2^{-n}$ for all $n \in {\mathbb N}$ and for $j=1,2$.
\end{description}
\end{proposition}

\begin{proof}
Let $\{v_n\}$ be as in Proposition~\ref{prop:1}. Then $\{v_n\}$ and
any subsequence thereof will satisfy (i). Upon passing to a
subsequence we can assume that
\begin{equation} \label{eq:2}
\|\lambda_{v_{m}}(S_i)
\lambda_{v_n}(S_j)-\lambda_{v_n}(S_j)\lambda_{v_{m}}(S_i)\| < 1/n,
\end{equation}
for all $m > n \ge 1$ and for all $i,j=1,2$. We claim that one can
find a sequence
$\{w_n\}$ of unitaries in $\F_2$ satisfying (ii) and (iii)
above---provided that we again pass to a subsequence of
$\{v_n\}$. It suffices to show that
for each $\delta > 0$ there exists a natural number $n$ such that for
each natural number $m > n$ there is a unitary $w \in \F_2$ for which
 $$\|w \lambda_{v_{m}}(S_j) w^* - \lambda_{v_n}(S_j)\| < \delta, \qquad
\|w S_j w^*-S_j\| < \delta$$ for $j=1,2$.
We give an indirect proof of the latter statement. If it were false, then there
 would exist $\delta >0$ and a sequence $1 \le n_1 < n_2 < n_3 <
 \cdots$ such that one of
$$\|w \lambda_{v_{n_{k+1}}}(S_i) w^* - \lambda_{v_{n_k}}(S_i)\|,
\qquad \|wS_iw^*-S_i\|,$$
$i=1,2$, is greater than $\delta$ for every $k$ and for all unitaries $w$ in
$\F_2$. We proceed
to show that this will lead to a contradiction.

Choose a free ultrafilter $\omega$ on $\mathbb N$ and consider the
relative commutant
$\O_2' \cap (\O_2)_\omega$ inside the ultrapower
$(\O_2)_\omega$. This $C^*$-algebra is purely infinite
and simple (see \cite[Proposition 3.4]{KirPhi:classI}).  Consider the
unital $*$-homomorphisms
$\eta_1,\eta_2 \colon
\O_2 \to \O_2' \cap (\O_2)_\omega$ given by
\begin{align*}
\eta_1(x) & = \pi_\omega\big(\lambda_{v_{n_2}}(x), \lambda_{v_{n_3}}(x),
\lambda_{v_{n_4}}(x), \dots\big), \\
\eta_2(x) & = \pi_\omega\big(\lambda_{v_{n_1}}(x), \lambda_{v_{n_2}}(x),
\lambda_{v_{n_3}}(x), \dots\big),
\end{align*}
$x \in \O_2$, where $\pi_\omega \colon \ell^\infty(\O_2) \to
(\O_2)_\omega$ is the quotient mapping. The images
of $\eta_1$ and $\eta_2$ commute by \eqref{eq:2}. Put
$$u = \eta_2(S_1)\eta_1(S_1)^*+\eta_2(S_2)\eta_1(S_2)^* =
\pi_\omega(v_{n_1}v_{n_2}^*,v_{n_2}v_{n_3}^*,v_{n_3}v_{n_4}^*,
\dots),$$
and notice that $u$ is a unitary element in $\O_2' \cap (\F_2)_\omega
\subseteq \O_2' \cap (\O_2)_\omega$. Use Corollary~\ref{cor:1.5} to
obtain a sequence $\{w_n\}$ of unitaries in $C^*(\eta_1(\F_2),u)
\subseteq \O_2' \cap (\F_2)_\omega$
such that $w_n\eta_1(S_j)w_n^* \to \eta_2(S_j)$ for $j=1,2$.
By \cite[Lemma 2.5]{KirRor:pi2} there is a single
unitary $w$ in $\O_2' \cap (\F_2)_\omega$ such that
$w\eta_1(S_j)w^* = \eta_2(S_j)$ for
$j=1,2$ (and hence such that $w\eta_1(x)w^* = \eta_2(x)$ for all $x \in
\O_2$).

Each unitary element in the ultrapower $(\F_2)_\omega$ lifts to a
unitary element in $\ell^\infty(\F_2)$, so we can write
$$w = \pi_\omega(w_1,w_2,w_3, \dots),$$
where each $w_n$ is a unitary element in $\F_2$. This establishes the
desired contradiction, as
$$\lim_{n \to \omega} \|S_jw_n - w_nS_j\| = 0, \qquad \lim_{n\to
  \omega} \|w_n \lambda_{v_{n_{k+1}}}(S_j) w_n^* -
\lambda_{v_{n_k}}(S_j) \| = 0,$$
for $j=1,2$ and for all $k$.
\end{proof}

\begin{theorem} \label{thm:O2}
There is a unitary element $u \in \F_2$ such that
the relative commutant  $\lambda_u(\O_2)' \cap \O_2$ contains a
unital copy of $\O_2$.
\end{theorem}

\begin{proof}
Let $\{v_n\}$ and $\{w_n\}$ be as in Proposition~\ref{prop:2} and
define endomorphisms on $\O_2$ by
$$\lambda_n(x) = w_1w_2 \cdots w_n \lambda_{v_{n+1}}(x)w_n^* \cdots
w_2^*w_1^*, \quad \rho_n(x) = w_1w_2 \cdots w_n  x w_n^* \cdots
w_2^*w_1^*,$$
for $x \in \O_2$. Then
$$\|\lambda_n(S_j)-\lambda_{n-1}(S_j)\| < 2^{-n}, \qquad
\|\rho_n(S_j)-\rho_{n-1}(S_j)\| < 2^{-n}$$
for $j=1,2$,
and $\lambda_n(x)\rho_n(y)-\rho_n(y)\lambda_n(x) \to 0$ for all $x,y
\in \O_2$. Using that
$$w \lambda_u(x) w^* = \lambda_{wu\varphi(w)^*}(x)$$
whenever $w$ is a unitary in $\O_2$ and $x \in \O_2$,
we see that $\lambda_n = \lambda_{u_n}$ for some
unitary $u_n$ in $\F_2$. It follows from the estimates above that
the sequences $\{\lambda_n(S_j)\}$ and $\{\rho_n(S_j)\}$,
$j=1,2$, and hence also the
  sequence $\{u_n\}$, are Cauchy and therefore convergent. Let
  $\lambda \colon \O_2 \to \O_2$ and $\rho \colon \O_2 \to \O_2$
  be the (pointwise-norm) limits of the sequences $\{\lambda_n\}$ and
  $\{\rho_n\}$, respectively, and let $u \in \F_2$ be the limit of
  the sequence $\{u_n\}$. Then $\lambda = \lambda_u$ and the images of
  $\lambda$ and $\rho$ commute.
\end{proof}

\begin{corollary} \label{cor:a}
There is a unitary $v \in \O_2$ such that $\lambda_v(\F_2) \subseteq
\F_2$ but $v \notin \F_2$.
\end{corollary}

\begin{proof} Let $u \in \O_2$ be as in Theorem~\ref{thm:O2} and take a unitary
  element $z$ in $\lambda_u(\O_2)' \cap \O_2$ that does not belong to
  $\F_2$. Put $v = zu$. Then $v$ does not belong to $\F_2$, and
  $\lambda_u$ and $\lambda_v$ coincide on $\F_2$ by
  Proposition~\ref{uv}, whence $\lambda_v$ maps $\F_2$ into itself.
\end{proof}

\begin{corollary} \label{cor:b} There is a unital $*$-homomorphism $\sigma
  \colon \O_2 \otimes \O_2 \to \O_2$ such that $\sigma(\F_2 \otimes
  \F_2) \subseteq \F_2$.
\end{corollary}

\begin{proof} Take $\lambda \colon \O_2 \to \O_2$ and $\rho \colon
  \O_2 \to \O_2$ as in the proof of Theorem~\ref{thm:O2}. Recall that
  $\lambda$ and $\rho$ have commuting images and that $\lambda(\F_2)
  \subseteq \F_2$ and  $\rho(\F_2) \subseteq \F_2$. We can therefore
  define a  $*$-homomorphism $\sigma \colon \O_2 \otimes \O_2 \to \O_2$
  by
$$\sigma(x \otimes y) = \lambda(x)\rho(y), \qquad x,y \in \O_2.$$
Then
$$\sigma(\F_2 \otimes \F_2) = \lambda(\F_2)\rho(\F_2) \subseteq \F_2.$$
\end{proof}

\noindent We know that $\O_2 \otimes \O_2$ and $\O_2$ are isomorphic,
but we do not know if one can find an \emph{isomorphism}  $\sigma
  \colon \O_2 \otimes \O_2 \to \O_2$ such that $\sigma(\F_2 \otimes
  \F_2)$ is contained in (or better, equal to) $\F_2$.

\section{Endomorphisms preserving the canonical $\mathbf{UHF}$-subalgebra}\label{epcs}

Below, $\phi$ denotes the standard left inverse of $\varphi$, i.e.,
the unital, completely positive map given by
$\phi(x):=\frac{1}{n}\sum S_i^* x S_i$, $x\in\O_n$.

\begin{lemma}\label{f1}
Let $u \in \U(\O_n)$, then the following conditions are equivalent:
\begin{description}
\item[\hspace{2mm}{\rm (i)}] $\phi(u) \in \U(\O_n)$;
\item[\hspace{1mm}{\rm (ii)}] $u \in \varphi(\O_n)$;
\item[{\rm (iii)}] $S_i^* u S_i = S_j^* u S_j \in \U(\O_n)$, for all $i,j \in \{1,\ldots,n\}$.
\end{description}
\end{lemma}
\begin{proof}
(i) $\Rightarrow$ (ii):
it follows from (i) that $u$ lies in the multiplicative domain of $\phi$ and therefore,
by Choi's theorem, $\phi(S_i u) = \phi(S_i) \phi(u)$, that is $u S_i = S_i \phi(u)$
for all $i=1,\ldots,n$. Thus, $u = \varphi(\phi(u))$.\\
The implications (ii) $\Rightarrow$ (iii) and (iii) $\Rightarrow$ (i) are obvious.
\end{proof}

\begin{lemma} \label{L1}
For $v,w\in\U(\O_n)$ the following three conditions are equivalent.
\begin{description}
\item[\hspace{2mm}{\rm (i)}] Endomorphisms $\lambda_v$ and $\lambda_{w}$ coincide on $\F_n$.
\item[\hspace{1mm}{\rm (ii)}] For each $k\geq1$ we have $w_k^* v_k \in \varphi^k(\O_n)$.
\item[{\rm (iii)}] There exists a sequence of unitaries $z_k\in\U(\O_n)$ such that
$z_1=\phi(w^*v)$ and $z_{k+1}=\phi(w^*z_k v)$ for all $k\geq1$.
\end{description}
\end{lemma}
\begin{proof}
The endomorphisms $\lambda_v$ and $\lambda_w$ coincide on $\F_n$ if and only if they coincide
on each $\F_n^k$. Now if $\alpha$ and $\beta$ are two multi-indices of length $k$ then
$\lambda_v(S_\alpha S^*_\beta)=v_k S_\alpha S^*_\beta v_k^*$ and $\lambda_w(S_\alpha S^*_\beta)
=w_k S_\alpha S^*_\beta w_k^*$. Thus $\lambda_v(S_\alpha S^*_\beta)=\lambda_w(S_\alpha S^*_\beta)$
for all such $\alpha,\beta$ if and only if $w_k^*v_k$ is in the commutant of $\F_n^k$, that is
when $w_k^* v_k \in \varphi^k(\O_n)$. Now it easily follows from Lemma \ref{f1} that this holds
for all $k$ if and only if condition (iii) above is satisfied.
\end{proof}

\begin{proposition} \label{P2}
If $w\in\U(\O_n)$ then $\lambda_w(\F_n) \subseteq \F_n$ if and only if $\lambda_w$ and
$\lambda_{\alpha_t(w)}$ coincide on $\F_n$ for all $t \in {\mathbb R}$. This in turn takes
place if and only if one can inductively define unitaries $z^{(k)}_t$, $k \geq 1$,
$t \in {\mathbb R}$ by
$$ \varphi(z^{(1)}_t) = w^* \alpha_t(w), \;\;\;
({\rm Ad}w\circ\varphi)(z^{(k+1)}_t) = z^{(k)}_t. $$
Moreover, in that case $t \mapsto z^{(1)}_t$ ($t \in {\mathbb R}$) is a unitary $\alpha$-cocycle
in $\lambda_w(\F_n)' \cap \O_n$. Finally, if $\lambda_w(\F_n) \subset \F_n$ then
$w \in \F_n$ if and only if $\lambda_w$ and $\lambda_{\alpha_t(w)}$ have the same range for
all $t\in{\mathbb R}$.
\end{proposition}
\begin{proof}
Given a unitary $w$ in $\O_n$ one has, by a direct computation,
$$\lambda_{\alpha_t(w)} = \alpha_t \circ \lambda_w \circ \alpha^{-1}_t \ , $$
for all $t \in {\mathbb R}$. Since $\F_n$ is precisely the fixed point algebra under
the gauge action, the first claim is now clear. The second equivalence in terms of
the existence of the unitaries $z_t^{(k)}$ is then deduced from Lemma \ref{L1} (see also
Remark \ref{oldcocycles}, below). Now
notice that if such unitaries exist one has, for any $s, t \in {\mathbb R}$,
\begin{align*}
\varphi(z_{t+s}^{(1)})
& = w^* \alpha_{t+s}(w) = w^* \alpha_t(\alpha_s(w)) = w^* \alpha_t(w w^* \alpha_s(w)) \\
& = w^* \alpha_t(w) \alpha_t(w^* \alpha_s(w)) =
\varphi(z_t^{(1)}) \alpha_t(\varphi(z_s^{(1)})) \ ,
\end{align*}
from which the cocycle equation for $z^{(1)}$ follows immediately, since
$\alpha$ and $\varphi$ commute. Moreover, $z^{(1)}_t \in \lambda_w(\F_n)' \cap \O_n$
for all $t \in {\mathbb R}$ by identity (\ref{rcUHF}).

Finally, suppose that $\lambda_w(\F_n) \subseteq \F_n$ and $\lambda_w(\O_n) =
\lambda_{\alpha_t(w)}(\O_n)$.
Then, for each $t$, define a map $\beta_t$ from $\O_n$ into itself via
$\lambda_{\alpha_t(w)}(x) = \lambda_w ( \beta_t (x) ), x \in \O_n$.
It must necessarily be that $\beta_t \in {\rm Aut}_{\F_n}(\O_n)$
and the argument in the proof of \cite[Proposition 2.1(b)]{Cun2} goes through.
\end{proof}

In particular, if $\lambda_w \in {\rm Aut}(\O_n)$ or, more generally,
$\lambda_w(\F_n)'\cap\O_n=\Complessi 1$ then $\lambda_w(\F_n) \subseteq \F_n$
if and only if $w \in \F_n$.

\begin{remark}\label{oldcocycles}
\rm Existence of unitaries $z^{(k)}_t$, as defined in Proposition \ref{P2}, is easily seen
to be equivalent to existence of unitaries $\tilde{z}^{(k)}_t$, $k \geq 1$, $t \in {\mathbb R}$,
defined inductively by
$$ \varphi(\tilde{z}^{(1)}_t) = w^* \alpha_t(w), \;\;\;
\varphi(\tilde{z}^{(k+1)}_t) = w^* \tilde{z}^{(k)}_t\alpha_t(w). $$
\end{remark}

\begin{proposition} \label{P3}
Let $w\in\U(\O_n)$ be such that $\lambda_w(\F_n^1) \subseteq \F_n$. Then the unitary
$\alpha$-cocycle $z^{(1)}_t := \phi(w^* \alpha_t(w))$ is a coboundary, i.e. there
exists a unitary $z$ such that $z^{(1)}_t = z \alpha_t(z^*)$ for all $t \in {\mathbb R}$.
\end{proposition}
\begin{proof}
Indeed, since $\lambda_w(\F_n^1) \subseteq \F_n$ there exists a unitary $u\in\F_n$ such that
$\lambda_w$ and $\lambda_u$ coincide on $\F_n^1$. In fact, we could take as $\lambda_u$ an inner
automorphism implemented by a unitary in $\F_n$. Then $w^* u$ commutes with $\F_n^1$, and thus
there exists a unitary $z$ such that $w^* u = \varphi(z)$. Now we have
$\varphi(z\alpha_t(z^*))=w^*u\alpha_t(u^*)\alpha_t(w)=w^*\alpha_t(w)$, since
$\alpha_t(u^*) = u^*$.
\end{proof}

\begin{proposition} \label{P4}
If $w$ is a unitary in $\O_n$ such that $w \D_n w^* \subseteq \F_n$ then $w^* \alpha_t(w) \in \D_n$
for all $t \in {\mathbb R}$. If, in addition, $\lambda_w(\F_n) \subseteq \F_n$
then $t \mapsto z^{(1)}_t$ is a one-parameter unitary group in $\lambda_w(\F_n)' \cap \D_n$.
\end{proposition}
\begin{proof}
By assumption, for any $x \in \D_n$ one has
$w x w^* = \alpha_t(w x w^*)$ for all $t \in {\mathbb R}$. Therefore, $\D_n$ being a MASA
in $\O_n$, $w^* \alpha_t(w) \in \D_n' \cap \O_n = \D_n$. Now, notice that $\D_n \cap
\varphi(\O_n) = \varphi(\D_n)$, so that indeed if $\lambda_w(\F_n) \subseteq \F_n$ the
cocycle given by Proposition \ref{P2} lies in $\D_n \subseteq \F_n$, and the conclusion
follows at once from the cocycle equation and Proposition \ref{P2}.
\end{proof}

Of course, the first part of the preceding proposition applies to all elements of the group $\S_n$,
as they normalize $\D_n$.

The following result is a slight reformulation of Proposition \ref{uv},
enhanced for our needs, put in a more symmetric form and
taking also into account Proposition \ref{Prc} and Lemma \ref{Laut}.

\begin{proposition} \label{P5}
Let $u$ and $w$ be two unitaries in $\O_n$.
If $\lambda_u$ and $\lambda_w$ coincide on $\F_n$ then, for every nonnegative integer $h$,
$$({\rm Ad} u \circ \varphi)^h (w u^*) \in \lambda_u(\F_n)' \cap \O_n \ . $$
Conversely, if $wu^* \in \lambda_u(\F_n)' \cap \O_n$ then $\lambda_u(x) = \lambda_w(x)$
for any $x \in \F_n$.
\end{proposition}
\begin{proof}
Concerning the first implication, by the above
it clearly suffices to show only the case $h = 0$.
%
%
Indeed, for every $k \geq 1$ one has
\begin{align*}
w u^* \lambda_u(S_{\alpha_1} \ldots S_{\alpha_k} S_{\beta_k}^* \cdots S_{\beta_1}^*) u w^*
& = w u^* u S_{\alpha_1} \ldots u S_{\alpha_k} S_{\beta_k}^* u^* \cdots S_{\beta_1}^* u^* u w^* \\
& = w  S_{\alpha_1}uS_{\alpha_2} \ldots u S_{\alpha_k} S_{\beta_k}^* u^* \cdots
S_{\beta_2}^* u^* S_{\beta_1}^* w^* \\
& = w S_{\alpha_1} \lambda_u(S_{\alpha_2} \ldots S_{\alpha_k} S_{\beta_k}^* \ldots S_{\beta_2}^*)
S_{\beta_1}^* w^* \\
& = w S_{\alpha_1} \lambda_w(S_{\alpha_2} \ldots S_{\alpha_k} S_{\beta_k}^* \ldots S_{\beta_2}^*)
S_{\beta_1}^* w^* \\
& = \lambda_w(S_{\alpha_1} \ldots S_{\alpha_k} S_{\beta_k}^* \cdots S_{\beta_1}^*) \\
& = \lambda_u(S_{\alpha_1} \ldots S_{\alpha_k} S_{\beta_k}^* \cdots S_{\beta_1}^*) \ .
\end{align*}
The opposite implication can be easily checked by induction on $k$,
just repeating the argument in Proposition \ref{uv}
after noticing that if $wu^* \in \lambda_u(\F_n)' \cap \O_n$
then, by Lemma \ref{Laut},
$wu^* = u\varphi(z)u^*$ for some unitary $z \in \lambda_u(\F_n)' \cap \O_n$,
that is $w = u \varphi(z)$.
\end{proof}

In particular, it follows
that if $w \in \U(\O_n) \setminus \F_n$ and there exists some $u \in \U(\F_n)$ such that
$\lambda_w$ and $\lambda_u$ coincide on $\F_n$ then $w$ must necessarily be of the form
$w = u \varphi(z)$ for some $z \in \lambda_u(\F_n)' \cap\O_n$,
which is exactly the situation discussed in section \ref{counterexample}.

\begin{corollary} \label{C6}
Let $w$ be a unitary in $\O_n$ and suppose that $\lambda_w(\F_n) = \F_n$.
Then $w \in \F_n$ and $\lambda_w \in {\rm Aut}(\O_n)$.
\end{corollary}
\begin{proof}
By Proposition \ref{P2}, $\lambda_w$ and $\lambda_{\alpha_t(w)}$ coincide on $\F_n$.
Therefore, by Proposition \ref{P5} one has
$w\alpha_t(w^*) \in \U(\F_n ' \cap \O_n) ={\mathbb T}$ and thus $w$ is an eigenvector for $\alpha$.
Hence $w$ belongs to $\F_n$ and the conclusion follows from \cite[Proposition 1.1 (a)]{CKS}.
\end{proof}

Combining \cite[Proposition 1.1 (a)]{CKS}, Proposition \ref{P2} and Corollary \ref{C6},
we obtain the following.

\begin{corollary} \label{C67}
For a unitary $w \in \O_n$, the following three conditions are equivalent:
\begin{description}
\item[\hspace{2mm}{\rm (i)}]  $\lambda_w(\F_n) = \F_n$;
\item[\hspace{1mm}{\rm (ii)}]  $\lambda_w \in {\rm Aut}(\O_n)$ and $w \in \F_n$;
\item[{\rm (iii)}] $\lambda_w \in {\rm Aut}(\O_n)$ and $\lambda_w(\F_n) \subseteq \F_n$.
\end{description}
\end{corollary}

\begin{corollary} \label{C7}
Assume that $\lambda_w$ is an endomorphism of $\O_n$ that restricts to the identity on $\F_n$.
Then $\lambda_w$ is a gauge automorphism.
\end{corollary}

\begin{proof}
By Corollary \ref{C6}, $\lambda_w \in {\rm Aut}_{\F_n}(\O_n) = \{\alpha_t : t \in {\mathbb R} \}$.
\end{proof}

\begin{corollary} \label{C8}
Let $w$ be a unitary in $\O_n$
such that $w^* \D_n w \subseteq \F_n$.
If $\lambda_w(\F_n) \subseteq \F_n$ then $w \alpha_t(w^*)
\in \lambda_w(\F_n)' \cap \D_n$ so that, in particular,
$w \in \F_n$ whenever $\lambda_w$ is irreducible in restriction to $\F_n$.
\end{corollary}

\begin{proof}
This readily follows from Propositions \ref{P2}, \ref{P4} and \ref{P5}.
\end{proof}

\begin{corollary}
Let $w \in \S_n$ be such that $\lambda_w(\D_n) = \D_n$ or, more generally, such that
$\D_n \subseteq \lambda_w(\F_n)$. Then $\lambda_w(\F_n) \subseteq \F_n$ if and only if $w \in \P_n$.
\end{corollary}

\begin{proof}
An element of $\S_n$ normalizes $\D_n$ and thus satisfies the first assumption
in the previous corollary.
Then the only nontrivial assertion follows from the fact that an endomorphism
$\lambda_w$ of $\O_n$ such that $\lambda_w(\F_n) \supseteq \D_n$ is necessarily irreducible
in restriction to $\F_n$ by an argument similar to the one in
\cite[Proposition 1.1]{CKS},
using the facts that $\D_n$ is a MASA in $\F_n$ and $\F_n$ is simple.
\end{proof}

\begin{example}\label{snpndn}
\rm In order to provide a simple example, we consider the following situation.
Let $w' \in \P_n$ be such that $\lambda_{w'}(\D_n) = \D_n$ but
$\lambda_{w'} \notin {\rm Aut}(\O_n)$ (many examples of such permutation unitaries
were provided in \cite{CS}). Let $w'' \in \S_n \setminus \P_n $ be such that
$\lambda_{w''} \in {\rm Aut}(\O_n)$ (e.g., an inner one), so that $\lambda_{w''}(\D_n) = \D_n$.
Set $\lambda_w := \lambda_{w'} \lambda_{w''}$, then $w = \lambda_{w'}(w'')w'
\in \S_n \setminus \P_n$ and $\lambda_w(\D_n) = \D_n$.
Such $\lambda_w$ is irreducible on $\O_n$ and $\lambda_w(\F_n) \not\subseteq \F_n$.
\end{example}

\begin{remark}\label{sndn}
\rm In view of the above, it would also be very useful to have a general criterion
for $w \in \S_n$ to satisfy
\begin{itemize}
\item[\hspace{2mm}{\rm (i)}] $\lambda_w(\D_n) = \D_n$, or
\item[\hspace{1mm}{\rm (ii)}] $\lambda_w \in {\rm Aut}(\O_n)$.
\end{itemize}
Such criteria for $w\in\P_n$ were given in \cite{CS,S}.
\end{remark}

\section{Analysis of cocycles for unitaries in $\S_n$}\label{sncocycles}

Let $w \in \S_n$. Then $w$ is of the form
$$w = \sum_{(\alpha,\beta)} S_\alpha S_\beta^* \ , $$
where the sum runs over a certain family $\J$ of pairs of multi-indices $(\alpha,\beta)$.
For convenience, we also introduce a set $\J_2 := \{\beta \ | \ (\alpha,\beta) \in \J\}$, which is
in bijective correspondence with $\J$ via $(\alpha,\beta) \leftrightarrow \beta$.
The fact that $w$ as above is unitary is equivalent to that both collections of the $P_\alpha$'s
and of the $P_\beta$'s form partitions of unity, i.e.
$$\sum_{(\alpha,\beta)} P_\alpha = \sum_{(\alpha,\beta)}P_\beta = 1 \ . $$
Then, for each $i=1,\ldots,n$, one has $S_i^* \sum P_\alpha S_i = 1$
and similarly for the $P_\beta$'s and therefore,
after summing over all $i$'s,
$$
\sum_{i=1}^n S_i^* \Big(\sum_{(\alpha,\beta)} P_\alpha\Big) S_i
= \sum_{i=1}^n S_i^* \Big(\sum_{(\alpha,\beta)} P_\beta\Big) S_i
= n 1 \ .
$$
Consequently, denoting by $\tilde\alpha$ (resp. $\tilde\beta$) the multi-index obtained from
$\alpha$ (resp. $\beta$) after deleting the first entry, we have
$$\sum_{(\alpha,\beta)} P_{\tilde\alpha} = \sum_{(\alpha,\beta)} P_{\tilde\beta} = n 1 \ . $$
In other words, both collections of projections $\{P_{\tilde\alpha}\}$ and $\{P_{\tilde\beta}\}$
form an $n$-covering of unity.

In the sequel, we repeatedly make use of Proposition \ref{P2} without further mention.
We compute for $t \in {\mathbb R}$
\begin{align*}
w^* \alpha_t (w)
& =
\Big(\sum_{(\alpha,\beta)} S_\beta S_\alpha^*\Big)
\Big(\sum_{(\alpha',\beta')} e^{it(|\alpha'|-|\beta'|)} S_{\alpha'}S_{\beta'}^*\Big) \\
& = \sum_{(\alpha,\beta)} e^{it(|\alpha|-|\beta|)} P_\beta
\in \U(\D_n)
\end{align*}
by orthogonality of the ranges of $S_\alpha$'s. Throughout the reminder of this section, we
assume that $\lambda_w(\F_n) \subseteq \F_n$. Therefore,
$w^* \alpha_t(w)$ must be a unitary in $\D_n \cap \varphi(\O_n) = \varphi(\D_n)$ and hence
$z^{(1)}_t := \phi(w^* \alpha_t (w))$ must be a unitary in $\D_n$.
We have
$$
z^{(1)}_t =
\frac{1}{n} \sum_{i=1}^n S_i^* \Big(\sum_{(\alpha,\beta)} e^{it(|\alpha|-|\beta|)} P_\beta \Big) S_i
= \frac{1}{n}\sum_{(\alpha,\beta)} e^{it(|\alpha|-|\beta|)} P_{\tilde\beta} \ .
$$
The last expression turns out to be unitary precisely when $|\alpha| - |\beta|$ is constant
over the classes of $\J_2$ with respect to the equivalence relation ``generated by nontrivial
overlaps of the $P_{\tilde\beta}$'s''. Namely, for $\beta, \beta' \in \J_2$, define
$$
\beta \sim \beta' \Leftrightarrow
\mbox{$\exists \beta_1 = \beta, \ldots, \beta_r = \beta' \in \J_2, \;\;
P_{\tilde\beta_s} P_{\tilde\beta_{s+1}} \neq 0$, $ \forall s = 1,\ldots,r-1$} \ .
$$
Thus, $z^{(1)}_t$ is unitary if and only if the function $\psi_1\colon\J_2\rightarrow{\mathbb Z}$ such
that $\psi_1(\beta)=|\alpha|-|\beta|$ is constant on the equivalence classes of relation $\sim$.
Unfortunately, such combinatorial analysis of "higher cocycles" $z^{(k)}_t$ quickly becomes
rather cumbersome. Thus from now on we make a simplifying assumption that for all $(\alpha,\beta)\in\J$
we have $|\alpha|-|\beta|\in\{-1,0,+1\}$.

\begin{proposition}\label{psifunctions}
Let $w = \displaystyle{\sum_{(\alpha,\beta)\in\J}} S_\alpha S_\beta^* \in \S_n$ be
such that $|\alpha|-|\beta|\in\{-1,0,+1\}$ for all $(\alpha,\beta)\in\J$. Then
$\lambda_w(\F_n) \subseteq \F_n$ if and only if there exists
a sequence of functions $\psi_k\colon \J_2 \to {\mathbb Z}$, $k=1,2,\ldots$, such that
\begin{itemize}
\item[(1)] $\psi_k$ is constant on the equivalence classes of relation $\sim$,
\item[(2)] $\psi_1(\beta) = |\alpha| - |\beta|$ and $\psi_{k+1}(\beta) = \psi_k(\beta')$,
where $(\alpha,\beta)\in\J$ and $\beta'$ is any element of $\J_2$ such that $\widetilde{\beta'}$ is
an initial segment of $\alpha$.
\end{itemize}
If such functions exist then
\begin{equation}\label{psiz}
z^{(k)}_t = \frac{1}{n} \sum_{(\alpha,\beta)\in\J} e^{it\psi_k(\beta)} P_{\tilde\beta}
\end{equation}
are unitary for all $k=1,2,\ldots$
\end{proposition}
\begin{proof}
By Proposition \ref{P2}, $\lambda_w(\F_n) \subseteq \F_n$ if and only if all "higher cocycles"
$z^{(k)}_t$, $k=1,2,\ldots$, are unitary. We show by induction on $k$ then under our hypothesis on
$w$ there exist functions $\psi_k$, $k=1,\ldots,m$ satisfying conditions (1) and (2) above if
and only if cocycles $z^{(k)}_t$, $k=1,\ldots,m$ are unitary and given by formula (\ref{psiz}).
Case $k=1$ is established just above this lemma. So assume the inductive hypothesis holds
for $k$. Then a direct calculation yields
\begin{equation}\label{zkp}
z^{(k+1)}_t = \phi(w^*z_t^{(k)}w) = \frac{1}{n} \sum_{(\alpha,\beta)\in\J}
e^{it\psi_k(\beta')} P_{\tilde\beta},
\end{equation}
where $\beta'$ is an element of $\J_2$ such that $\widetilde{\beta'}$ is an initial segment
of $\alpha$. Note that for another such element $\beta''$ we have $\psi_k(\beta')=
\psi_k(\beta'')$, since function $\psi_k$ is constant on equivalence classes of relation $\sim$.
Thus we can define $\psi_{k+1}(\beta)=\psi_k(\beta')$. Formula (\ref{zkp}) yields a unitary
if and only if function $\psi_{k+1}$ is constant on equivalence classes of $\sim$. This ends
the proof of the inductive step and the lemma.
\end{proof}

The conditions of Proposition \ref{psifunctions} can be given the following graphical
interpretation. Let $w = \sum_{(\alpha,\beta)\in\J} S_\alpha S_\beta^* \in \S_n$ be
such that $|\alpha|-|\beta|\in\{-1,0,+1\}$ for all $(\alpha,\beta)\in\J$. We associate with $w$ a
finite directed graph $E_w$ as follows. Vertices of $E_w$ are the equivalence classes of
relation $\sim$. Given two vertices $a_1$, $a_2$, there is a single edge from $a_1$ to $a_2$
if and only if there exist $(\alpha,\beta)\in\J$ with $\beta$ in the equivalence class $a_1$
and $\beta'\in\J_2$ in the equivalence class $a_2$ such that $\widetilde{\beta'}$ is an
initial segment of $\alpha$. We denote by $E_w^k$ the collection of all directed paths in
$E_w$ of length $k$, and by $E_w^k(a)$ the collection of those such paths which begin at
vertex $a$.

If the function $\psi_1$ (corresponding to $w$) is constant on the equivalence classes of
$\sim$ then we can assign labels from $\{-1,0,+1\}$ to vertices of $E_w$ in such a way
that the label of $a$ is $\psi_1(\beta)$ for $\beta$ in the equivalence class $a$. Now the
remaining conditions of Lemma \ref{psifunctions} are equivalent to the following
{\em path condition}:

\vspace{3mm}\noindent
For each vertex $a$ and for each $k\in{\mathbb N}$ the ranges of all directed paths in
$E_w^k(a)$ have the same labels.

\vspace{3mm}
Since the graph $E_w$ is finite, we obtain the following:

\begin{corollary}\label{finiteness}
Let $w = \sum_{(\alpha,\beta)\in\J} S_\alpha S_\beta^* \in \S_n$ be
such that $|\alpha|-|\beta|\in\{-1,0,+1\}$ for all $(\alpha,\beta)\in\J$.
Then there exists $r\in{\mathbb N}$ such that
$$ \lambda_w(\F_n^r)\subseteq\F_n \; \Rightarrow \; \lambda_w(\F_n)\subseteq\F_n. $$
\end{corollary}

\begin{example}\label{firstexample}
\rm Let $u\in\P_2^4$ and $v\in\S_2$ be as given in Example \ref{countercuntz}.
Set $w=vu$. Then the corresponding graph $E_w$ looks as follows.

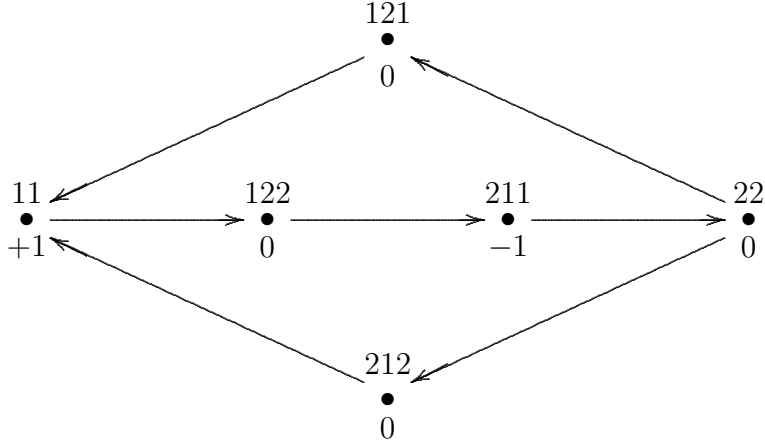
\[ \beginpicture
\setcoordinatesystem units <1.6cm,1.2cm>
\setplotarea x from -1 to 7, y from -2.5 to 2.5

\put {$\bullet$} at 0 0
\put {$\bullet$} at 2 0
\put {$\bullet$} at 4 0
\put {$\bullet$} at 6 0
\put {$\bullet$} at 3 2
\put {$\bullet$} at 3 -2

\setlinear

\plot 0.2 0  1.8 0 /
\plot 2.2 0  3.8 0 /
\plot 4.2 0  5.8 0 /
\plot 0.2 0.2  2.8 1.8 /
\plot 0.2 -0.2  2.8 -1.8 /
\plot 3.2 1.8  5.8 0.2 /
\plot 3.2 -1.8  5.8 -0.2 /

\arrow <0.235cm> [0.2,0.6] from 1.6 0 to 1.8 0
\arrow <0.235cm> [0.2,0.6] from 3.6 0 to 3.8 0
\arrow <0.235cm> [0.2,0.6] from 5.6 0 to 5.8 0
\arrow <0.235cm> [0.2,0.6] from 0.5 -0.4 to 0.2 -0.2
\arrow <0.235cm> [0.2,0.6] from 0.5 0.4 to 0.2 0.2
\arrow <0.235cm> [0.2,0.6] from 3.5 -1.6 to 3.2 -1.8
\arrow <0.235cm> [0.2,0.6] from 3.5 1.6 to 3.2 1.8

\put{$11$} at 0 0.3
\put{$122$} at 2 0.3
\put{$211$} at 4 0.3
\put{$22$} at 6 0.3
\put{$121$} at 3 2.3
\put{$212$} at 3 -1.6

\put{$+1$} at 0 -0.3
\put{$0$} at 2 -0.3
\put{$-1$} at 4 -0.3
\put{$0$} at 6 -0.3
\put{$0$} at 3 1.6
\put{$0$} at 3 -2.3

\endpicture \]
\end{example}

We denote by $A_w$ the smallest $({\rm Ad}w\circ\varphi)$-invariant $C^*$-subalgebra
of $\D_n$ that contains $\{z_t^{(1)}:t\in\Reali\}$.

\begin{proposition}\label{invsub}
Let $w = \sum_{(\alpha,\beta)\in\J} S_\alpha S_\beta^* \in \S_n$ be
such that $|\alpha|-|\beta|\in\{-1,0,+1\}$ for all $(\alpha,\beta)\in\J$.
Then $\lambda_w(\F_n)\subseteq\F_n$ if and only if $z_t^{(1)}$ is a unitary
cocycle and the algebra $A_w$ is finite dimensional. In that case, $A_w$ is the
$C^*$-algebra generated by all cocycles $\{z_t^{(k)}:k\in\Naturali\}$.
\end{proposition}
\begin{proof}
If $A_w$ is finite dimensional then ${\rm Ad}w\circ\varphi$ is its automorphism.
If, in addition, $z_t^{(1)}$ is unitary then this immediately implies existence of
unitary cocycles $z_t^{(k)}$ for all $k\in\Naturali$.

Conversely, if $\lambda_w(\F_n)\subseteq\F_n$ then the $C^*$-algebra generated by
all cocycles $\{z_t^{(k)}:k\in\Naturali\}$ is finite dimensional, since it is contained
in $C^*(\{P_{\tilde{\beta}}:\beta\in J_2\})$. It follows that ${\rm Ad}w\circ\varphi$
is an automorphism of this algebra.
\end{proof}



\newpage\noindent
Roberto Conti \\
Department of Mathematics \\
University of Rome 2 Tor Vergata\\
Via della Ricerca Scientifica, 00133 Rome, Italy \\
E-mail: conti@mat.uniroma2.it \\

\smallskip \noindent
Mikael R{\o}rdam \\
Department of Mathematical Sciences \\
The University of Copenhagen \\
Universitetspark 5, DK--2100 Copenhagen, Denmark \\
E-mail: rordam@math.ku.dk \\

\smallskip \noindent
Wojciech Szyma{\'n}ski\\
Department of Mathematics and Computer Science \\
The University of Southern Denmark \\
Campusvej 55, DK-5230 Odense M, Denmark \\
E-mail: szymanski@imada.sdu.dk


\begin{thebibliography}{99}

\bibitem{BraKisRorSto:shift}
O.~Bratteli, A.~Kishimoto, M.~R{\o}rdam, and E.~St{\o}rmer, \emph{{The crossed
  product of a UHF-algebra by a shift}}, {Ergod. Th. \& Dynam. Sys.}
  \textbf{13} (1993), 615--626.

\bibitem{CKS} R. Conti, J. Kimberley and W. Szyma{\'n}ski,
{\em More localized automorphisms of the Cuntz algebras},
arXiv:0808.2843, to appear in Proc. Edinburgh Math. Soc.

\bibitem{CP} R. Conti and C. Pinzari,
{\it Remarks on the index of endomorphisms of Cuntz algebras},
J. Funct. Anal. {\bf 142} (1996), 369--405.

\bibitem{CS} R. Conti and W. Szyma{\'n}ski,
{\em Labeled trees and localized automorphisms of the Cuntz algebras},
arXiv:0805.4654.

\bibitem{Cun1} J. Cuntz, {\it Simple $C^*$-algebras generated by isometries},
Commun. Math. Phys. {\bf 57} (1977), 173--185.

\bibitem{Cun2} J. Cuntz, {\it Automorphisms of certain simple $C^*$-algebras},
in {\it Quantum fields-algebras-processes}, ed. L. Streit, Springer 1980.

\bibitem{KirPhi:classI}
E.~Kirchberg and N.~C. Phillips, \emph{{Embedding of exact $C^*$-algebras into
  ${\mathcal{O}}_2$}}, J. reine angew. Math. \textbf{525} (2000), 17--53.

\bibitem{KirRor:pi2}
E.~Kirchberg and M.~R{\o}rdam, \emph{{Infinite non-simple $C^*$-algebras:
  absorbing the Cuntz algebra $\mathcal{O}_\infty$}}, Adv. Math.
  \textbf{167} (2002), 195--264.

\bibitem{L} M. Laca, {\it Endomorphisms of ${\mathcal B}({\mathcal H})$ and Cuntz
algebras}, J. Operator Theory {\bf 30} (1993), 85--108.

\bibitem{Ne} V. Nekrashevych, {\it Cuntz-Pimsner algebras of group actions},
J. Operator Theory {\bf 52} (2004), 223--249.

\bibitem{Ror:cuntz}
M. R{\o}rdam, {\it Classification of inductive limits of Cuntz algebras},
J. reine angew. Math. {\bf 440} (1993), 175--200.

\bibitem{Ror:o2o2}
M. R{\o}rdam,
{\em A short proof of Elliott's theorem: $\O_2\otimes\O_2\cong\O_2$},
C. R. Math. Rep. Acad. Sci. Canada {\bf 16} (1994), 31--36.

\bibitem{S} W. Szyma{\'n}ski,
{\it On localized automorphisms of the Cuntz
algebras which preserve the diagonal subalgebra}, in `New Development of
Operator Algebras', R.I.M.S. K\^{o}ky\^{u}roku {\bf 1587} (2008), 109--115.

\bibitem{Y} D. Yang, {\it Endomorphisms and modular theory of 2-graph $C^*$-algebras},
arXiv:0907.1129, to appear in Indiana Univ. Math. J.

\end{thebibliography}
\end{document}